\documentclass[12pt]{elsarticle}

\usepackage{graphicx}
\usepackage{float}
\usepackage{subfigure}
\usepackage{amsfonts}
\usepackage{amsmath}
\usepackage{amsthm}
\usepackage{amssymb}
\usepackage{bm}
\usepackage{color}
\usepackage[T1]{fontenc}
\usepackage[utf8]{inputenc}
\usepackage{comment}
\usepackage{mathtools}
\usepackage{xcolor}

\usepackage{lineno}
\usepackage{natbib}

	{\left\lbrace \begin{array}{@{} l @{} }}%
	{ \end{array}\right.}

\newtheorem{pro}{Proposition}
\newtheorem{cor}{Corollary}

\newtheorem{thm}{Theorem}

\newtheorem{remark}{Remark}
\newtheorem{conj}{Conjecture}

\newcommand{\bbN}{\mathbb{N}}

\newcommand{\calT}{\mathcal{T}}
\newcommand{\calU}{\mathcal{U}}
\newcommand{\calX}{\mathcal{X}}
\newcommand{\de}{\,\mathrm{d}}

\newcommand{\bs}{\boldsymbol}

\definecolor{orange}{RGB}{255,127,0}

\begin{document}

\begin{frontmatter}

\title{On $(\beta,\gamma)$-Chebyshev functions and points \\of the interval}

\author{Stefano De Marchi$^{*,\dagger}$}
\ead{demarchi@math.unipd.it, corresponding author}

\author{Giacomo Elefante$^{**}$}
\ead{giacomo.elefante@unifr.ch}

\author{Francesco Marchetti$^{*}$}
\ead{francesco.marchetti@math.unipd.it}

\address{$^{*}$Dipartimento di Matematica \lq\lq Tullio Levi-Civita\rq\rq, Universit\`a di Padova, Italy; \\$^{**}$Département de Mathématiques, Université de Fribourg, Switzerland;}

\begin{abstract}
In this paper, we introduce the class of $(\beta,\gamma)$-Chebyshev functions and corresponding points, which can be seen as a family of {\it generalized} Chebyshev polynomials and points. For the $(\beta,\gamma)$-Chebyshev functions, we prove that they are orthogonal in certain subintervals of $[-1,1]$ with respect to a weighted arc-cosine measure. In particular we investigate the cases where they become polynomials, deriving new results concerning classical Chebyshev polynomials of first kind. Besides, we show that subsets of Chebyshev and Chebyshev-Lobatto points are instances of $(\beta,\gamma)$-Chebyshev points. We also study the behavior of the Lebesgue constants of the polynomial interpolant at these points  on varying the parameters $\beta$ and $\gamma$.
\end{abstract}


\end{frontmatter}


\section{Introduction}
Chebyshev polynomials have been long-investigated in scientific literature and they have been considered in various fields, e.g. in function approximation \cite{Shuman11}, partial differential equations \cite{Haidvogel79}, cryptography \cite{Bergamo05}, distributed consensus \cite{Montijano13}, group theory \cite{Bircan12}, cosmography \cite{Cappozziello18} as well as optimal control problems \cite{KDK12}. Different types of Chebyshev polynomials have been studied and many related properties have been reported (for a complete overview the interested reader may refer to \cite{Mason02,Rivlin74}).

In the recent literature, Chebyshev polynomials still represent a prolific research topic. For example, although generalizations of such polynomials have been already proposed for example in \cite{Laine80,Oliveira73}, more recent ones in \cite{Borzov19,Cesarano14,Hassani19}. Pseudo-Chebyshev functions of rational degree $p/q$ have been also recently studied in \cite{Caratelli20, Cesarano19}. As in the standard setting with integer degrees, it has been proved that the family pseudo-Chebyshev functions satisfies a recurrence relation and solves a certain differential equation, similar to the classical Chebyshev polynomials, however it retains an orthogonality property in an interval of the real axis if and only if $q=2$. Furthermore, new identities about Chebyshev polynomials have been derived also in \cite{Boussayoud19, Zhang18}. 

The zeros of Chebyshev polynomials, the Chebyshev points, are of large interest in literature. They represent a preferable choice for interpolation due to their well conditioning and fast convergence (cf. e.g. \cite{Cheney98,Rivlin03}). In fact, Chebyshev points retain a logarithmic growth for the norm of the interpolant operator. i.e. the Lebesgue constant \cite{Brutman78,Ehlich66}. Furthermore, for their good properties, they are widely-adopted, for example, in numerical quadrature \cite{MasjedJamei05} or in the solution of differential equations \cite{Sweilam16}.  

We fix some notations. Let $\Omega=[-1,1]$ and $n\in\mathbb{N}$. The Chebyshev polynomials of the first kind $\{T_n\}_{n=0,1,\dots}$ are defined as \begin{equation*}
    T_n(x)=\cos(n\arccos{x}),\; x\in\Omega \,.
\end{equation*}
$T_n$ is indeed an algebraic polynomial of degree $n$ thanks to the  change of variable $x=\cos(t)$ and the Viète formulae for the cosine. They are a family of orthogonal polynomials on $\Omega$ with respect to the weight function $w(x)=(1-x^2)^{-1/2}$. The zeros of $T_n$, namely the {\it Chebyshev points} (of the first kind), is the set
\begin{equation*}
    \mathcal{T}_n=\bigg\{\cos\bigg(\frac{(2j-1)\pi}{n}\bigg)\bigg\}_{j=1,\dots,n},\quad n\in\mathbb{N},
\end{equation*}
which are all reals and inside $\Omega$. 

Let $\mathcal{X}_n\coloneqq \{x_0,\dots,x_n\}$ be a set of distinct points in $\Omega$ and $\mathrm{L}\coloneqq\{\ell_0,\dots,\ell_n\}$ be the Lagrange polynomials 
\begin{equation*}
\ell_i(x)\coloneqq \prod_{\substack{j=0 \\ j\neq i}}^n {\frac{x-x_j}{x_i-x_j}},\; i=0,\dots,n,\;x\in\Omega,
\end{equation*}
the Lebesgue function is
\begin{equation*}
    \lambda(\mathcal{X}_n;x)=\sum_{i=0}^n{|\ell_i(x)|},\;x\in\Omega,
\end{equation*}
and its maximum over $\Omega$ is the corresponding {\it Lebesgue constant}
\begin{equation*}
    \Lambda(\mathcal{X}_n,\Omega)=\max_{x\in \Omega}\lambda(\mathcal{X}_n;x),
\end{equation*}
which is an indicator both of the conditioning and the stability of the interpolation process. Note that the Lebesgue constant depends only on the choice of the interpolation nodes, and therefore many efforts have been made in finding sets of nodes whose Lebesgue constant grows \textit{slowly} with $n$; we refer e.g. to \cite{Cheney98,Rivlin03}. 

Furthermore, also the set of \textit{Chebyshev-Lobatto} (CL) points 
\begin{equation*}
    \mathcal{U}_{n+1}=\bigg\{\cos\bigg(\frac{j\pi}{n}\bigg)\bigg\}_{j=0,\dots,n},
\end{equation*}
which consists of the zeros of the polynomial
\begin{equation*}
 \overline{T}_{n+1}(x)=\frac{(1-x^2)}{n}\frac{\partial }{\partial x}T_n(x),\quad x\in\Omega,   
\end{equation*}
is widely-adopted being as well $\Lambda(\mathcal{U}_{n+1},\Omega)=\mathcal{O}(\log{n})$ \cite{McCabe73}. We refer to the recent survey \cite{Ibrahimoglu16} for further details concerning the Lebesgue constant of Chebyshev, CL and various other sets of points.

In this work, we introduce a new family of functions in $\Omega$, namely the $(\beta,\gamma)$\textit{-Chebyshev} functions, which can be considered as a \textit{generalization} of Chebyshev polynomials. Indeed the family includes the classical Chebyshev polynomials as a particular case. 

After investigating such new functions and providing various theoretical results, we drive our attention to the corresponding sets of $(\beta,\gamma)$-Chebyshev and $(\beta,\gamma)$-CL points, analyzing the Lebesgue constant of these points from a theoretical point of view by verifying the results through extensive numerical experiments. We point out that such points can be characterized as \textit{mapped} equispaced points, and so they can be studied in the framework of the recently proposed \textit{fake nodes} (cf. \cite{Berrut20,DeMarchi20b,DeMarchi20,DeMarchi21}).

The paper layout and our main contributions follows.
\begin{itemize}
    \item 
    In Section \ref{sez_betta}, we introduce what we call the $(\beta,\gamma)$-Chebyshev functions and related points. In particular, in Theorem \ref{thm_orto} we prove that such functions are orthogonal in a subinterval of $\Omega$. Moreover, we analyze for which choice of the parameters $\beta,\gamma$ they reduce to polynomials. In doing so, we show how subsets of classical Chebyshev and CL points are included in our general framework. Furthermore, Corollary \ref{cor_orto} presents a new result concerning the orthogonality of standard Chebyshev polynomials of the first kind whose degree is a multiple of a fixed natural number.
    \item
    Section \ref{sez_mapped} is devoted to show how the new sets of nodes can be obtained by mapping sets of equispaced points through the Kosloff Tal-Ezer map (cf. \cite{Adcock16, Kosloff93}). Notice that in Proposition \ref{prop_equiv} we use this characterization to link $(\beta,\gamma)$-Chebyshev and $(\beta,\gamma)$-CL points.
    \item
    In Section \ref{sez_leb}, we investigate on the behavior of the Lebesgue constant of the interpolant at these points. More precisely, we show how the parameters $\beta$ and $\gamma$ influence the growth of the Lebesgue constant with respect to the degree $n$.
    \item
    Finally, in Section \ref{sez_conclusions} we draw some final considerations and discuss further developments.
\end{itemize}

\section{The $(\beta,\gamma)$-Chebyshev functions and related zeros}\label{sez_betta}
\subsection{The general case}
On $\Omega$, let us consider the $(\beta,\gamma)$-Chebyshev functions (of the first kind) defined as
\begin{equation}\label{eq:betagammacheb}
    T_n^{\beta,\gamma}(x)\coloneqq \cos\bigg(\frac{2n}{2-\beta-\gamma}\bigg(\arccos{x}-\frac{\gamma\pi}{2}\bigg)\bigg),\;x\in\Omega,
\end{equation}
where $\beta,\gamma\in[0,2)$, $\beta+\gamma<2$, $n\in\mathbb{N}$. We point out that in general $T_n^{\beta,\gamma}$ is not a polynomial and $T_n^{0,0}=T_n$ is the classical Chebyshev polynomials of the first kind. 

\begin{figure}[H]
  \centering
  \includegraphics[width=0.49\linewidth]{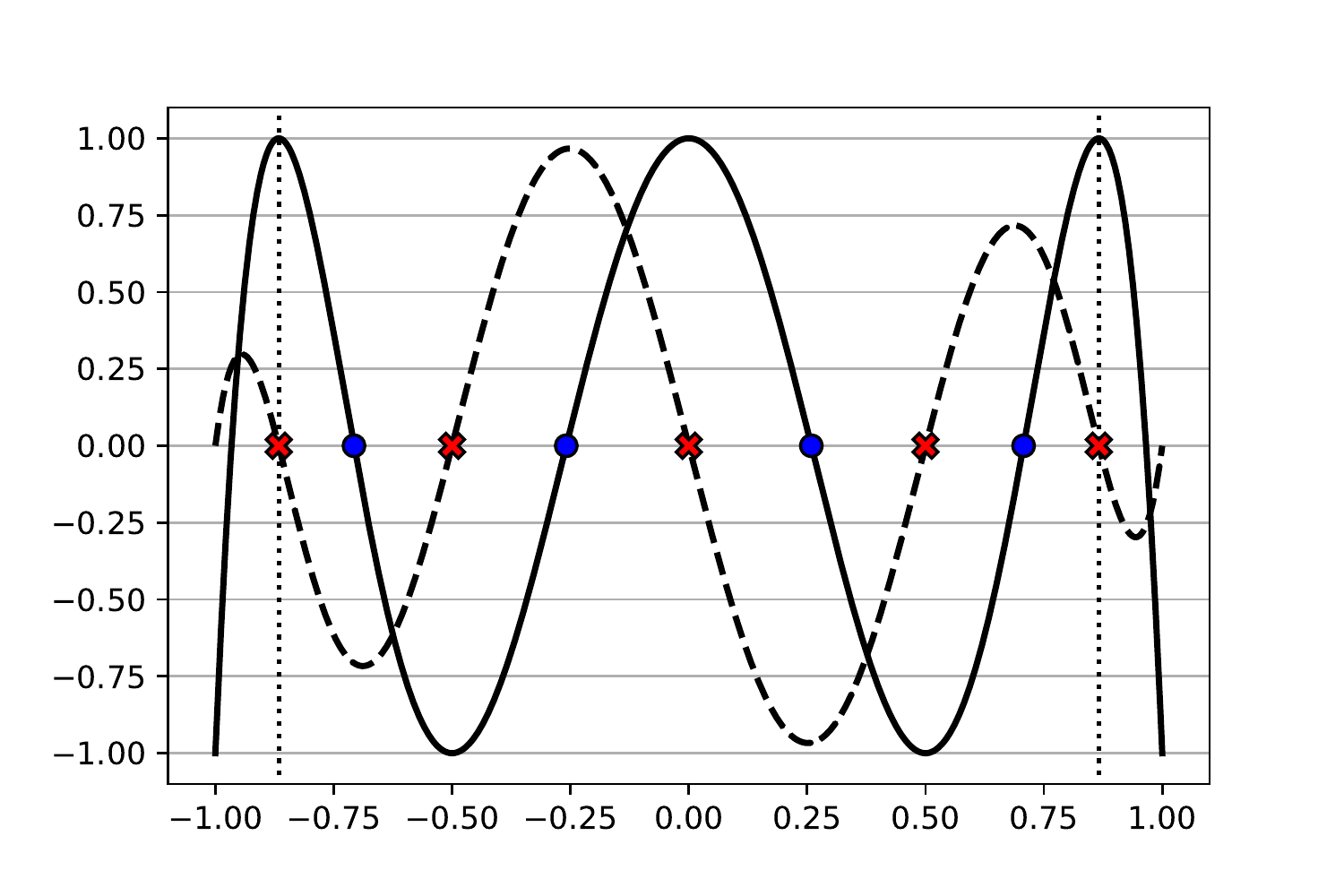}  
  \includegraphics[width=0.49\linewidth]{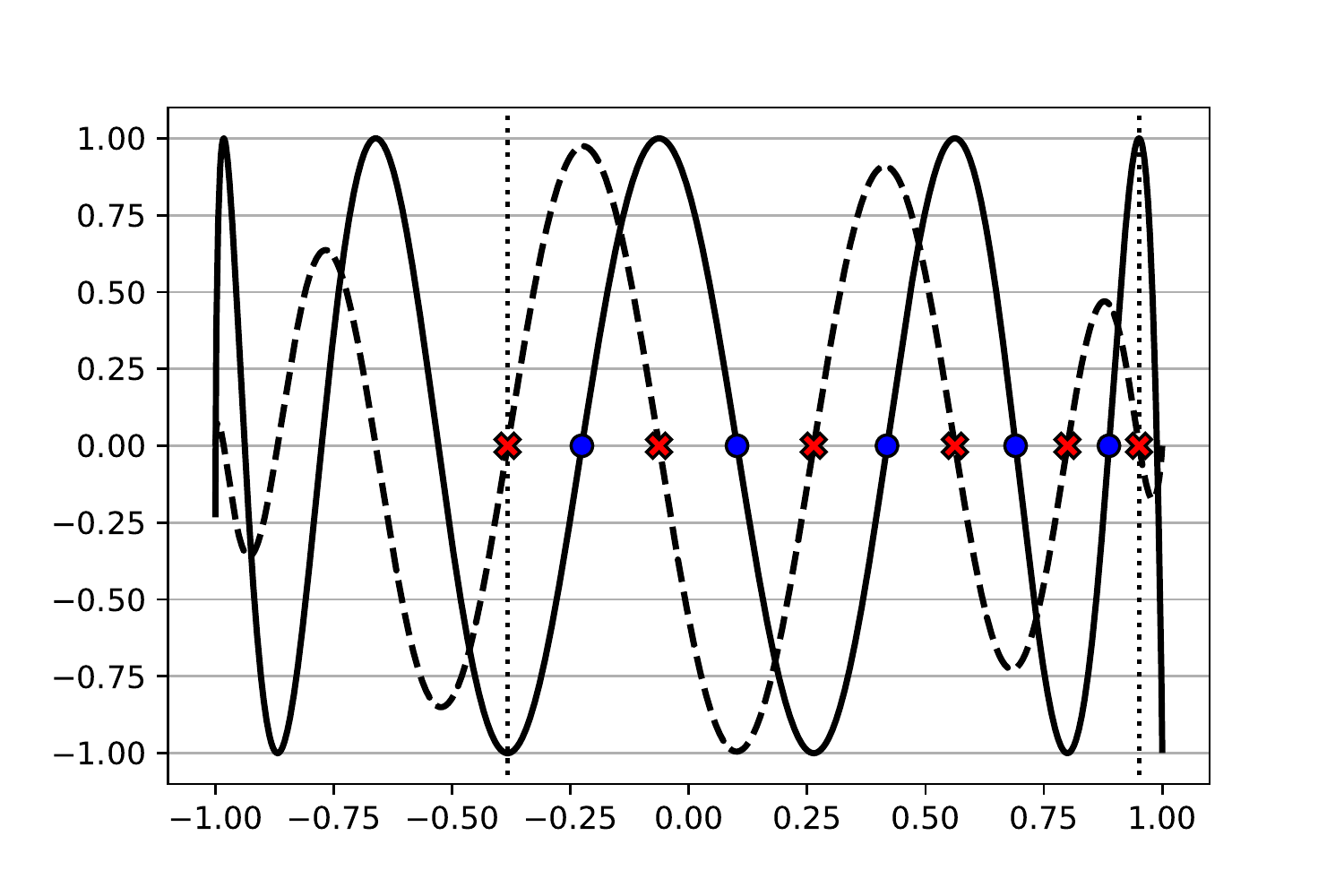}  
\caption{The functions $T_n^{\beta,\gamma}$ (solid line) and $\overline{T}_{n+1}^{\beta,\gamma}$ (dashed line), the sets $\mathcal{T}^{\beta,\gamma}_{n}$ (blue circles) and $\mathcal{U}^{\beta,\gamma}_{n+1}$ (red crosses), the set $\Omega_{\beta,\gamma}$ delimited by dotted vertical lines. Left: $n=4$, $\beta=\gamma=1/3$. Right: $n=5$, $\beta=3/4$, $\gamma=1/5$.}
\label{fig1}
\end{figure}

The set of zeros of the function $T_n^{\beta,\gamma}$ in $\Omega_{\beta,\gamma}\coloneqq[-\cos(\beta\pi/2),\cos(\gamma\pi/2)]\subseteq \Omega$, that is
\begin{equation*}
    \mathcal{T}^{\beta,\gamma}_n\coloneqq\bigg\{\cos\bigg(\frac{(2-\beta-\gamma)(2j-1)\pi}{4n}+\frac{\gamma\pi}{2}\bigg)\bigg\}_{j=1,\dots,n},
\end{equation*}
is what we call the \textit{$(\beta,\gamma)$-Chebyshev points}.
Moreover, the extrema points of $T_n^{\beta,\gamma}$ in $\Omega_{\beta,\gamma}$ are
\begin{equation*}
    \cos\bigg(\frac{(2-\beta-\gamma)j\pi}{2n}+\frac{\gamma\pi}{2}\bigg),\;j=1,\dots,n-1.
\end{equation*}
Similarly, we can define the set of \textit{$(\beta,\gamma)$-Chebyshev-Lobatto} ($(\beta,\gamma)$-CL) points as
\begin{equation*}
    \mathcal{U}^{\beta,\gamma}_{n+1}\coloneqq\bigg\{\cos\bigg(\frac{(2-\beta-\gamma)j\pi}{2n}+\frac{\gamma\pi}{2}\bigg)\bigg\}_{j=0,\dots,n}.
\end{equation*}
We note that the elements of $\mathcal{U}^{\beta,\gamma}_{n+1}$ are zeros of the function 
\begin{equation*}
    \overline{T}_{n+1}^{\beta,\gamma}(x)= \frac{2-\beta-\gamma}{2n} (1-x^2)\frac{\partial }{\partial x}T_n^{\beta,\gamma}(x),\;x\in\Omega.
\end{equation*}
More precisely, $\mathcal{U}^{\beta,\gamma}_{n+1}$ coincides with the set of zeros of $\overline{T}_{n+1}^{\beta,\gamma}$ if and only if $\beta=\gamma=0$. In Figure \ref{fig1}, we display the functions and corresponding points for some values of $n,\beta,\gamma$.
Now, we prove a result concerning a symmetric property of $(\beta,\gamma)$-Chebyshev functions.
\begin{pro}\label{prop:symmetry}
Let $n\in\mathbb{N}_{>0}$ and $\nu\in[0,2[$. Then, for $x\in\Omega$,
\begin{equation*}
    T_n^{\nu,0}(x)= (-1)^n\:T_n^{0,\nu}(-x),
\end{equation*}
and
\begin{equation*}
    \overline{T}_{n+1}^{\nu,0}(x)= (-1)^n\:\overline{T}_{n+1}^{0,\nu}(-x).
\end{equation*}
\end{pro}
\begin{proof}
By using the identity $\arccos{(-x)}=\pi-\arccos{x}$ and the addition formula for the cosine, we get
\begin{equation*}
    \begin{split}
        T_n^{0,\nu}(-x) &= \cos\bigg(\frac{2n}{2-\nu}\bigg(\pi-\arccos{x}-\frac{\nu\pi}{2}\bigg)\bigg)\\
        &= \cos\bigg(n\pi-\frac{2n}{2-\nu}\arccos{x}\bigg)\\
        &= \cos(n\pi)\cos\bigg(\frac{2n}{2-\nu}\arccos{x}\bigg)\\
        &= (-1)^n T_n^{\nu,0}(x),
    \end{split}
\end{equation*}
because of the fact that $\sin(n\pi)=0$ for any $n\in\mathbb{N}_{>0}$.\\
Moreover, for $x\in\Omega$,
\begin{align*}
    \overline{T}_{n+1}^{\nu,0}(x)&= \frac{2-\beta-\gamma}{2n} (1-x^2)\frac{\partial }{\partial x}T_n^{\nu,0}(x) \\
    &=(-1)^n \frac{2-\beta-\gamma}{2n} (1-x^2)\frac{\partial }{\partial x}T_n^{0,\nu}(-x)= (-1)^n \overline{T}_{n+1}^{0,\nu}(-x).
\end{align*}
This concludes the proof.
\end{proof}
\begin{cor}\label{cor:symmetry}
    In the hypotheses of Proposition \ref{prop:symmetry}, we have
    \begin{align*}
        \bar{x}\in \mathcal{T}^{\nu,0}_{n} \quad \text{if and only if}\quad -\bar{x}\in \mathcal{T}^{0,\nu}_{n}, \\
        \bar{x}\in \mathcal{U}^{\nu,0}_{n} \quad \text{if and only if}\quad  -\bar{x}\in \mathcal{U}^{0,\nu}_{n}.
    \end{align*}
\end{cor}
These properties are a direct consequence of the results in Proposition \ref{prop:symmetry}.

\begin{figure}[H]
  \centering
  \includegraphics[width=0.49\linewidth]{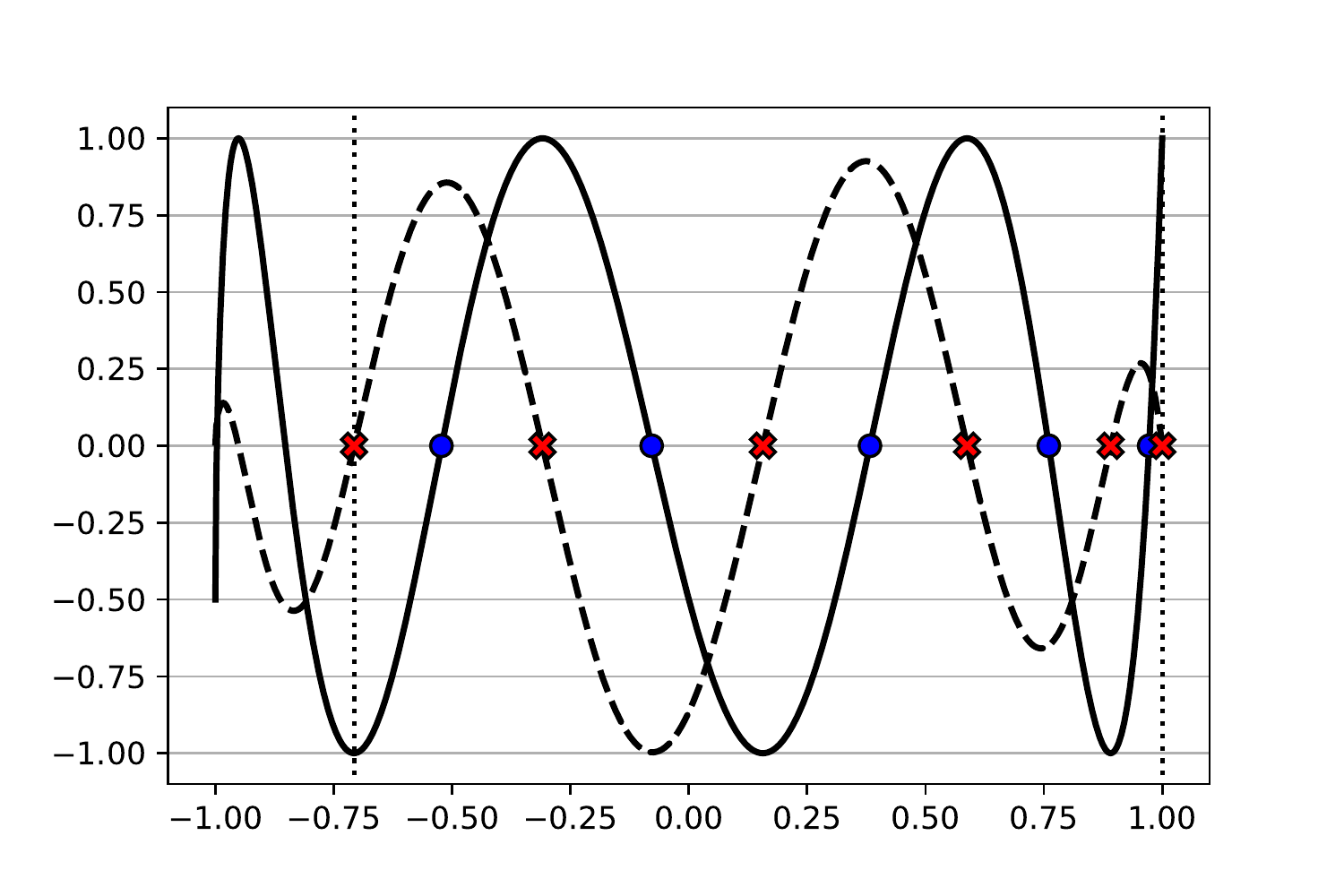}  
  \includegraphics[width=0.49\linewidth]{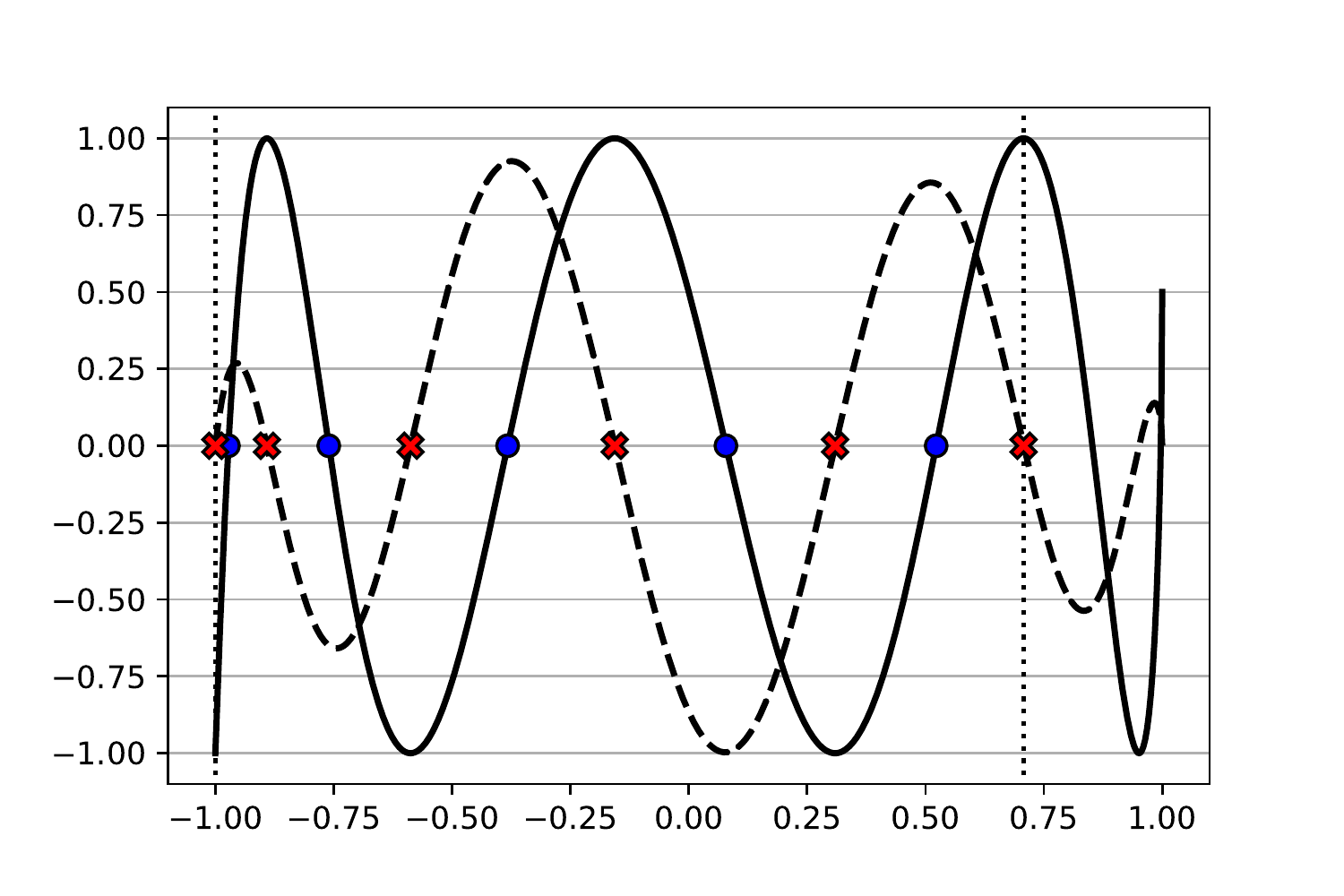}  
\caption{The functions $T_n^{\beta,\gamma}$ (solid line) and $\overline{T}_{n+1}^{\beta,\gamma}$ (dashed line), the sets $\mathcal{T}^{\beta,\gamma}_{n}$ (blue circles) and $\mathcal{U}^{\beta,\gamma}_{n+1}$ (red crosses), the set $\Omega_{\beta,\gamma}$ delimited by dotted vertical lines. Left: $n=5$, $\beta=1/2$, $\gamma=0$. Right: $n=5$, $\beta=0$, $\gamma=1/2$.}
\label{fig2}
\end{figure}

\begin{figure}[H]
  \centering
  \includegraphics[width=0.49\linewidth]{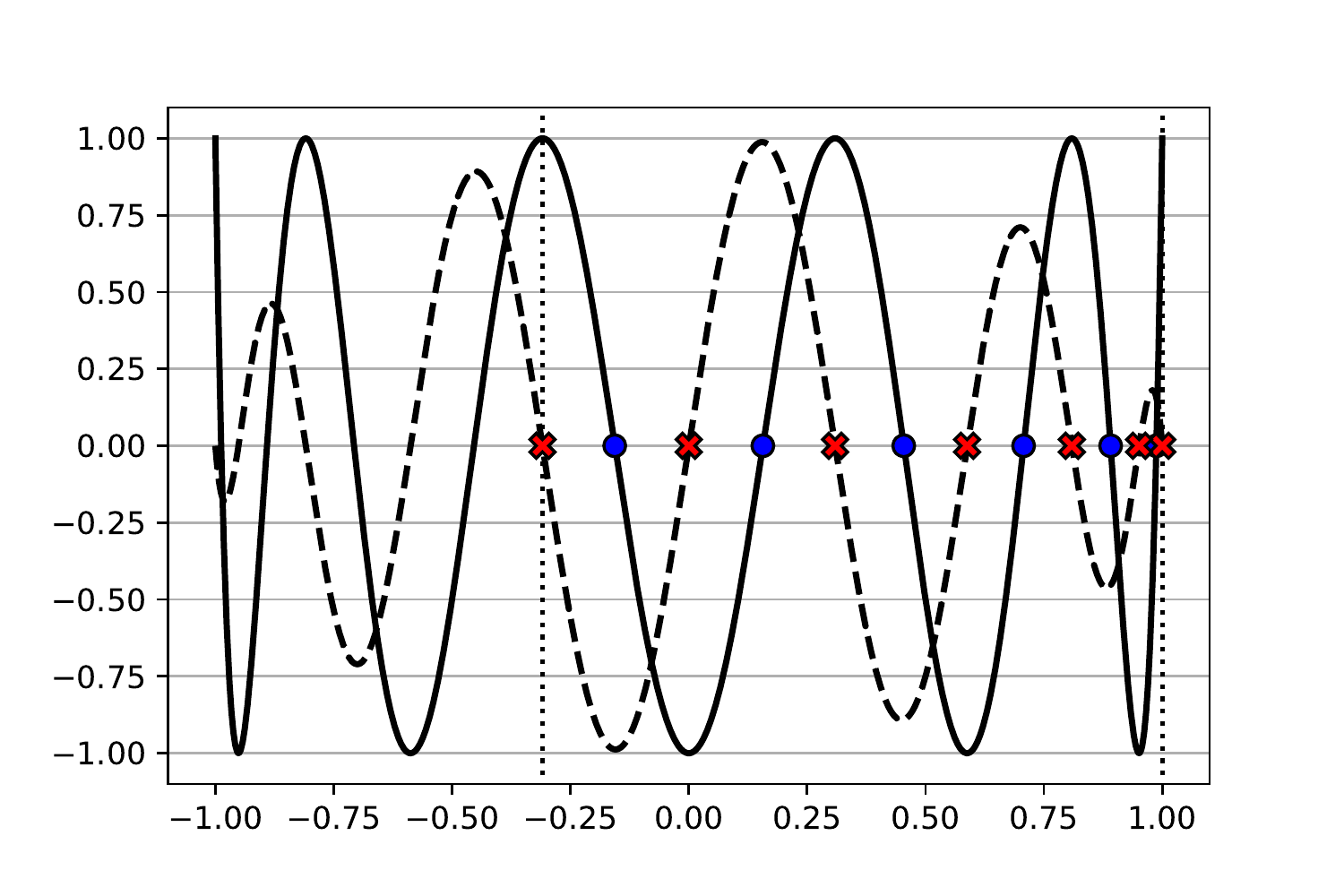}  
  \includegraphics[width=0.49\linewidth]{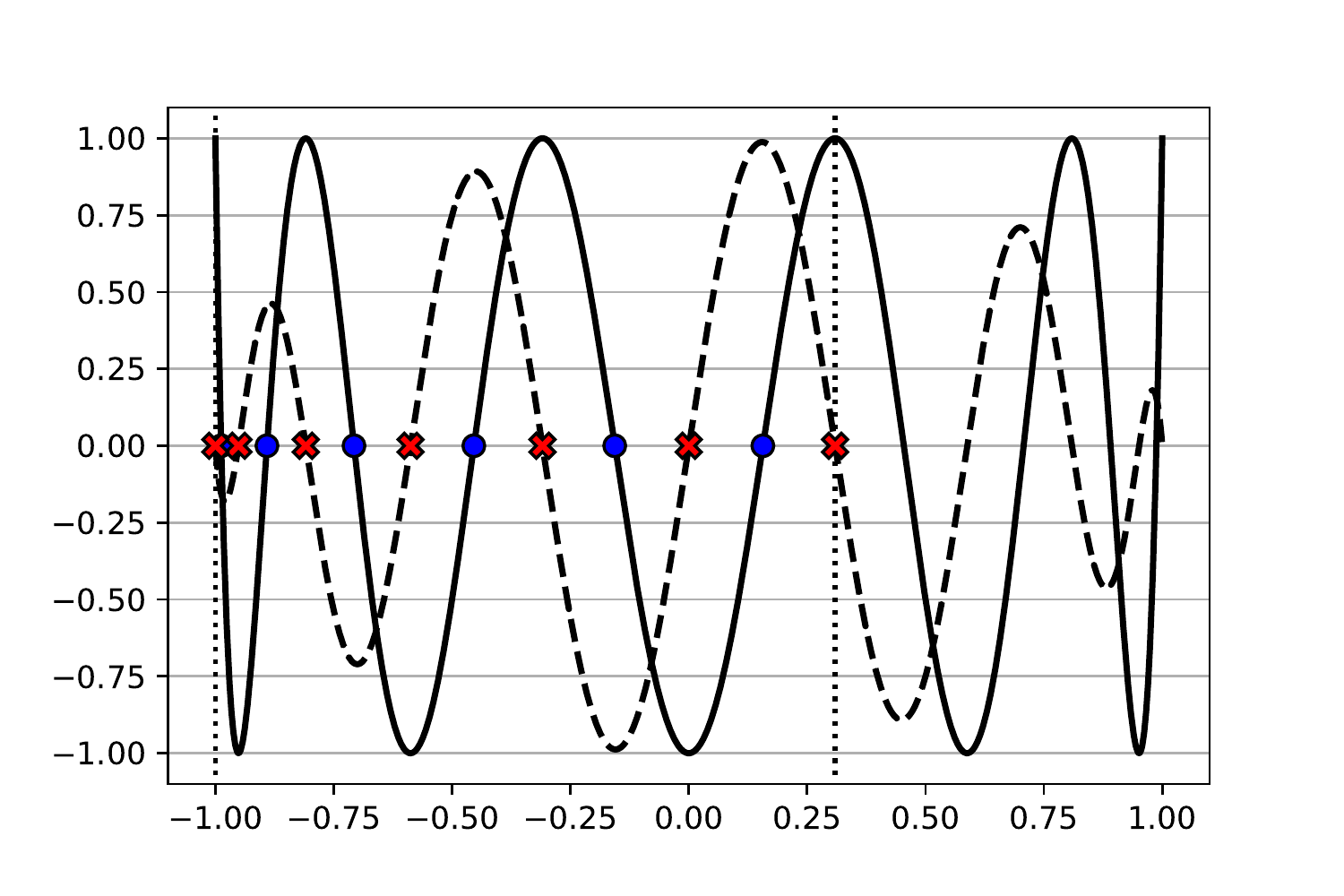}  
\caption{The functions $T_n^{\beta,\gamma}$ (solid line) and $\overline{T}_{n+1}^{\beta,\gamma}$ (dashed line), the sets $\mathcal{T}^{\beta,\gamma}_{n}$ (blue circles) and $\mathcal{U}^{\beta,\gamma}_{n+1}$ (red crosses), the set $\Omega_{\beta,\gamma}$ delimited by dotted vertical lines. Left: $n=6$, $\beta=4/5$, $\gamma=0$. Right: $n=6$, $\beta=4/5$, $\gamma=0$.}
\label{fig3}
\end{figure}

In Figures \ref{fig2} and \ref{fig3} we show the above symmetric properties for some values of $n,\beta,\gamma$.
This family of functions satisfies a recurrence formula.
\begin{pro}
The functions $\{T_n^{\beta,\gamma}\}_{n=0,1,\dots}$ satisfy
\begin{equation*}
\begin{split}
    & T_{0}^{\beta,\gamma}(x)=1,\quad T_{1}^{\beta,\gamma}(x)= \cos\bigg(\frac{2}{2-\beta-\gamma}\bigg(\arccos{x}-\frac{\gamma\pi}{2}\bigg)\bigg),\\ & T_{n+1}^{\beta,\gamma}(x)+T_{n-1}^{\beta,\gamma}(x)=2T_{1}^{\beta,\gamma}(x)T_{n}^{\beta,\gamma}(x)
    \end{split}
\end{equation*}
for $x\in\Omega$.
\end{pro}
\begin{proof}
Letting $\theta=\arccos{x}-\frac{\gamma\pi}{2}$, $x\in\Omega$, by using the addition formulae for the cosine we get
\begin{equation*}
    \cos\bigg(\frac{2(n+1)\theta}{2-\beta-\gamma}\bigg)+\cos\bigg(\frac{2(n-1)\theta}{2-\beta-\gamma}\bigg)=2\cos\bigg(\frac{2\theta}{2-\beta-\gamma}\bigg)\cos\bigg(\frac{2n\theta}{2-\beta-\gamma}\bigg),
\end{equation*}
which concludes the proof.
\end{proof}

Furthermore, we have the following orthogonality result.
\begin{thm}\label{thm_orto}
    The functions $\{T_n^{\beta,\gamma}\}_{n=0,1,\dots}$ are orthogonal on $\Omega_{\beta,\gamma}$ with respect to the weight function
    \begin{equation*}
    w^{\beta,\gamma}(x)=\frac{2}{(2-\beta-\gamma)\sqrt{1-x^2}},\; x\in \Omega_{\beta,\gamma},
    \end{equation*} 
    having
    \begin{equation*}
        \int_{\Omega_{\beta,\gamma}}{T_r^{\beta,\gamma}(x)T_s^{\beta,\gamma}(x)w^{\beta,\gamma}(x)\de x}=
        \begin{dcases} 0 & \textrm{if } r\neq s,\\
        \pi & \textrm{if }r=s=0,\\
        \frac{\pi}{2} & \textrm{if }r=s\neq 0.\end{dcases}
    \end{equation*}
\end{thm}
\begin{proof}
As well-known,
    \begin{equation*}
        \int_{0}^{\pi}{\cos(r\theta)\cos(s\theta)\de \theta}=
        \begin{dcases} 0 & \textrm{if } r\neq s,\\
        \pi & \textrm{if }r=s=0,\\
        \frac{\pi}{2} & \textrm{if }r=s\neq 0.\end{dcases}
    \end{equation*}
Then, the result follows by the change of variable $x=\cos((2-\beta)\theta/2+\gamma\pi/2)$ that maps $\Omega$ into $\Omega_{\beta,\gamma}$.
\end{proof}

\subsection{The polynomial case}
In what follows, we analyze under which choices of the parameters the function $ T_n^{\beta,\gamma}$ is a polynomial. We remark that thanks to Proposition \ref{prop:symmetry}, we do not need to discuss the case $\beta=0,\;\gamma>0$.\\ 
\subsubsection{Case $\beta>0,\;\gamma=0$}
Consider $\beta=p/q,\; p,q\in\mathbb{N}_{>0}$ and independent of $n$. From \eqref{eq:betagammacheb}, we then require
\begin{equation*}
    \frac{2}{2-\frac{p}{q}}\in\mathbb{N}_{>0}\Longleftrightarrow (2q-p)m=2q
\end{equation*}
for a given $m\in \mathbb{N}_{>1}$.
Hence, we obtain
\begin{equation}\label{eq:betam}
    p = \frac{2q(m-1)}{m} \quad \text{ and }\quad \beta_m \coloneqq \frac{2(m-1)}{m} = 2-\frac{2}{m}.
\end{equation}
Therefore, for $x\in\Omega$
\begin{equation*}
    T_n^{\beta_m,0}(x)=\cos\bigg(\frac{2n}{2-\beta_m}\arccos{x}\bigg)=\cos(mn\arccos{x})=T_{mn}(x),
\end{equation*}
which implies
\begin{equation*}
    \begin{split}
        & \mathcal{T}^{\beta_m,0}_n=\bigg\{\cos\bigg(\frac{(2j-1)\pi}{2mn}\bigg)\bigg\}_{j=1,\dots,n},\\
        & \mathcal{U}^{\beta_m,0}_{n+1}=\bigg\{\cos\bigg(\frac{j\pi}{mn}\bigg)\bigg\}_{j=0,\dots,n}.
    \end{split}
\end{equation*}
We observe that $\beta_m\in[1,2[$ for every $m\in \mathbb{N}_{>1}$. From Theorem \ref{thm_orto} we get the following corollary which shows that the Chebyshev polynomials (of the first kind) with degree and weight function that are a multiple of a fixed $m\in\bbN_{>1}$, satisfy this additional orthogonality property.

\begin{cor}\label{cor_orto}
    Let $m\in \mathbb{N}_{>1}$, the polynomials $\big\{T_n^{\beta_m,0}\big\}_{n=0,1,\dots}=\{T_{mn}\}_{n=0,1,\dots}$ are orthogonal on $\Omega_{\beta_m,0}$ with respect to the weight function
    \begin{equation*}
    w^{\beta_m,0}(x)=\frac{m}{\sqrt{1-x^2}},\;x\in\Omega_{\beta_m,0}.
    \end{equation*} 
\end{cor}
\begin{proof}
The proof directly follows from Theorem \ref{thm_orto}.
\end{proof}

Let us now consider the case $\beta$ depending linearly on $n$. In particular,
\begin{equation*}
    \frac{2n}{2-\frac{p}{q}}\in\mathbb{N}_{>0}\Longleftrightarrow (2q-p)m=2qn,\;m\in \mathbb{N}_{>0},
\end{equation*}
that is, 
\begin{equation}\label{eq_beta_poly}
    p = \frac{2q(m-n)}{m}  \quad \text{ and }\quad  \beta_{m,n} \coloneqq \frac{2(m-n)}{m} = 2-\frac{2n}{m},
\end{equation}
where we assume $n<m$ so that $\beta_{m,n}\in]0,2[$.
Then, for a fixed $m\in\mathbb{N}_{>1}$ and $n=0,\dots,m-1$ the corresponding functions and points are
\begin{equation*}
    T_n^{\beta_{m,n},0}(x)=\cos\bigg(\frac{2n}{2-\beta_{m,n}}\arccos{x}\bigg)=\cos(m\arccos{x})=T_{m}(x),
\end{equation*}
and
\begin{equation}\label{eq_punti_rimanenti}
    \begin{split}
        & \mathcal{T}^{\beta_{m,n},0}_n=\bigg\{\cos\bigg(\frac{(2j-1)\pi}{2m}\bigg)\bigg\}_{j=1,\dots,n},\\
        & \mathcal{U}^{\beta_{m,n},0}_{n+1}=\bigg\{\cos\bigg(\frac{j\pi}{m}\bigg)\bigg\}_{j=0,\dots,n}.
    \end{split}
\end{equation}
Where $\mathcal{T}^{\beta_{m,n},0}_n$ and $\mathcal{U}^{\beta_{m,n},0}_{n+1}$ are subsets of classical Chebyshev and CL points (see Theorem \ref{thm_sottoinsiemi} below). The case $m=n+1$ is of particular interest for us and is discussed in Section \ref{sez_leb}.

\subsubsection{Case $\beta>0,\;\gamma=\beta$}\label{sec_beta=gamma}
In view of \eqref{eq:betam}, we take
\begin{equation}\label{eq:betam2}
    \beta=\gamma=\frac{\beta_m}{2}.
\end{equation}

\begin{pro}\label{prop_equilibrata}
Let $m,n\in\mathbb{N}_{>1}$ and let $\beta=\gamma=\frac{\beta_m}{2}$.
\begin{enumerate}
    \item 
    If $n\in 4\mathbb{N}_{>0}$, then $T_n^{\frac{\beta_m}{2},\frac{\beta_m}{2}}=T_{mn}$.
    \item
    If $n\in 2\mathbb{N}_{>0}\setminus 4\mathbb{N}_{>0}$, then
    \begin{equation*}
        T_n^{\frac{\beta_m}{2},\frac{\beta_m}{2}} = \begin{dcases} T_{mn} & \textrm{if } m-1\in 2\mathbb{N}_{>0},\\
        -T_{mn} & \textrm{if } m\in 2\mathbb{N}_{>0}.\end{dcases}
    \end{equation*}
    \item
    If $n$ is odd, then
    \begin{equation*}
        T_n^{\frac{\beta_m}{2},\frac{\beta_m}{2}} = \begin{dcases} T_{mn} & \textrm{if } m-1\in 4\mathbb{N}_{>0},\\
        -T_{mn} & \textrm{if } m-1\in 2\mathbb{N}_{>0}\setminus 4\mathbb{N}_{>0},\\
        -\sin(mn\arccos{x}) & \textrm{if } (m-1)n+1\in 4\mathbb{N}_{>0},\\
        \sin(mn\arccos{x}) & \textrm{if } (m-1)n-1\in 4\mathbb{N}_{>0}.        
        \end{dcases}
    \end{equation*}
\end{enumerate}
\end{pro}
\begin{proof}
For $x\in\Omega$ we have
\begin{equation*}
    \begin{split}
        T_n^{\frac{\beta_m}{2},\frac{\beta_m}{2}}(x)&= \cos\bigg(mn\bigg(\arccos{x}-\frac{\beta_m}{2}\frac{\pi}{2}\bigg)\bigg)\\
        & =\cos\bigg(mn\arccos{x}-\frac{(m-1)n\pi}{2}\bigg)\\
        & = \cos(mn\arccos{x})\cos\bigg(\frac{(m-1)n\pi}{2}\bigg)+\\
        & + \sin(mn\arccos{x})\sin\bigg(\frac{(m-1)n\pi}{2}\bigg).
    \end{split}
\end{equation*}
Then the three cases follow by evaluating the sine and cosine for the corresponding values of $m,n\in\bbN_{>1}$.
\end{proof}

We note that the zeros of $T_n^{\frac{\beta_m}{2},\frac{\beta_m}{2}}$ in $\Omega$ are, therefore, $mn+1$ Chebyshev-Lobatto points in the case where $n$ is odd and $(m-1)n\pm 1\in 4\bbN_{>0}$, otherwise they are $mn$ Chebyshev points.

Concerning $\left(\frac{\beta_m}{2},\frac{\beta_m}{2}\right)$-Chebyshev and $\left(\frac{\beta_m}{2},\frac{\beta_m}{2} \right)$-CL points, we obtain
\begin{equation*}
    \begin{split}
        & \mathcal{T}^{\frac{\beta_m}{2},\frac{\beta_m}{2}}_n=\bigg\{\cos\bigg(\frac{(2j-1)\pi}{2mn}+\frac{(m-1)\pi}{2m}\bigg)\bigg\}_{j=1,\dots,n},\\
        & \mathcal{U}^{\frac{\beta_m}{2},\frac{\beta_m}{2}}_{n+1}=\bigg\{\cos\bigg(\frac{j\pi}{mn}+\frac{(m-1)\pi}{2m}\bigg)\bigg\}_{j=0,\dots,n}.
    \end{split}
\end{equation*}
\begin{cor}
    Let $m\in\mathbb{N}_{>1}$ be an odd number. The polynomials
    \begin{equation*}
        \big\{T_n^{\frac{\beta_m}{2},\frac{\beta_m}{2}}\big\}_{n=0,1,\dots} = \begin{dcases} \{T_{mn}\}_{n=0,1,\dots} & \textrm{if } m-1\in 4\mathbb{N}_{>0},\\
        \{(-1)^n T_{mn}\}_{n=0,1,\dots} & \textrm{if } m-1\in 2\mathbb{N}_{>0}\setminus 4\mathbb{N}_{>0},
        \end{dcases}
    \end{equation*}
     are orthogonal in $\Omega_{\frac{\beta_m}{2},\frac{\beta_m}{2}}$ with respect to the weight function
    \begin{equation*}
    w^{\frac{\beta_m}{2},\frac{\beta_m}{2}}(x)=\frac{m}{\sqrt{1-x^2}},\;x\in\Omega_{\frac{\beta_m}{2},\frac{\beta_m}{2}}.
    \end{equation*} 
\end{cor}
\begin{proof}
The proof directly follows from Theorem \ref{thm_orto} and Proposition \ref{prop_equilibrata}.
\end{proof}

Finally, in view of \eqref{eq_beta_poly}, we consider
\begin{equation*}
    \beta = \gamma = \frac{\beta_{m,n}}{2},\;n<m.
\end{equation*}
\begin{pro}\label{prop_equilibrata_bis}
Let $m\in\mathbb{N}_{>1}$ be fixed and let $\beta = \gamma = \frac{\beta_{m,n}}{2}$ with $n<m$.
\begin{enumerate}
    \item 
    If $m-n\in 4\mathbb{N}_{>0}$, then $T_n^{\frac{\beta_{m,n}}{2},\frac{\beta_{m,n}}{2}}=T_{m}$.
    \item
    If $m-n\in 2\mathbb{N}_{>0}\setminus 4\mathbb{N}_{>0}$, then
    $T_n^{\frac{\beta_{m,n}}{2},\frac{\beta_{m,n}}{2}}=-T_{m}$.
    \item
    If $m-n+1\in 4\mathbb{N}_{>0}$, then
    $T_n^{\frac{\beta_{m,n}}{2},\frac{\beta_{m,n}}{2}}(x)=-\sin(m\arccos{x})$ for $x\in\Omega$.
    \item
    If $m-n-1\in 4\mathbb{N}_{>0}$, then
    $T_n^{\frac{\beta_{m,n}}{2},\frac{\beta_{m,n}}{2}}(x)=\sin(m\arccos{x})$ for $x\in\Omega$.
\end{enumerate}
\end{pro}
\begin{proof}
The proof is similar to that of Proposition \ref{prop_equilibrata}.
\end{proof}

\subsubsection{General case}
Let $\beta=1-p_1/q_1,\;\gamma=1-p_2/q_2,\;p_1,p_2\in\mathbb{N},\;q_1,q_2\in\mathbb{N}_{>0}$. In general, $ T_n^{\beta,\gamma}$ is a polynomial if
$$  \frac{2n}{\frac{p_1}{q_1}+\frac{p_2}{q_2}}\in\mathbb{N}_{>0},$$
that is
\begin{equation*}
    p_1q_2+p_2q_1\:|\:2q_1q_2n \quad \Longleftrightarrow \quad (p_1q_2+p_2q_1)m=2q_1q_2n,\;m\in \mathbb{N}_{>0}.
\end{equation*}
It is worthwhile to point out a particular choice of $\beta$ and $\gamma$ for which the $(\beta,\gamma)$-Chebyshev (Lobatto) points result in subsets of Chebyshev (Lobatto) points. 
\begin{thm}\label{thm_sottoinsiemi}
Let $\mathcal{T}_n=\{t_j\}_{j=1,\dots,n}$ and $\mathcal{U}_{n+1}=\{u_j\}_{j=0,\dots,n}$ be the set of Chebyshev and CL points respectively. Moreover, let $\kappa_1,\kappa_2\in\bbN$, $\bs{\kappa}\coloneqq (\kappa_1,\kappa_2),$ and $\beta_{\bs{\kappa}}=\frac{2\kappa_1}{n+\kappa_1+\kappa_2}$, $\gamma_{\bs{\kappa}}=\frac{2\kappa_2}{n+\kappa_1+\kappa_2}$. Then,
\begin{align*}
    \calT^{\beta_{\bs{\kappa}},\gamma_{\bs{\kappa}}}_{n} &= \calT_{n+\kappa_1+\kappa_2} \setminus \{ t_1,\dots,t_{\kappa_2-1}, t_{n+\kappa_2+1}, \dots, t_{n+\kappa_1+\kappa_2}\},\\
    \calU^{\beta_{\bs{\kappa}},\gamma_{\bs{\kappa}}}_{n+1} &=\calU_{n+\kappa_1+\kappa_2+1}\setminus \{ u_0,\dots,u_{\kappa_2-1}, u_{n+\kappa_2+1}, \dots, u_{n+\kappa_1+\kappa_2}\}.
\end{align*}
\end{thm}

\begin{proof}
We present the proof only for $(\beta,\gamma)$-CL points because for the $(\beta,\gamma)$-Chebyshev points is similar.

\begin{align*}
    \calU^{\beta_{\bs{\kappa}},\gamma_{\bs{\kappa}}}_{n+1} &= \cos\left( \frac{2-\beta_{\bs{\kappa}}-\gamma_{\bs{\kappa}}}{2n}j\pi + \frac{\gamma_{\bs{\kappa}} \pi}{2}\right), \quad j=0,\dots,n, \\
    &= \cos\left( \frac{ \frac{2n}{n+\kappa_1+\kappa_2}}{2n}j\pi + \frac{ \kappa_2 \pi}{n+\kappa_1+\kappa_2} \right), \\
    &= \cos\left( \frac{(j+\kappa_2)\pi}{n+\kappa_1+\kappa_2} \right), \\
    &=  \cos\left( \frac{l\,\pi}{n+\kappa_1+\kappa_2} \right), \quad l=\kappa_2,\dots,n+\kappa_2.
\end{align*}
\end{proof}

With the notation of Theorem \ref{thm_sottoinsiemi} and recalling the analysis carried out in Section \ref{sec_beta=gamma}, we highlight 
\begin{equation*}
    \beta=\gamma=\frac{\beta_{n+2,n}}{2}=1-\frac{n}{n+2}=\frac{2}{n+2}=\frac{2\kappa_1}{n+\kappa_1+\kappa_2},
\end{equation*}
with $\kappa_1=\kappa_2=1$ and in this case, we obtain indeed the sets
\begin{equation}\label{eq_puntirimasti2}
\begin{split}
    & \calU^{\frac{\beta_{n+2,n}}{2},\frac{\beta_{n+2,n}}{2}}_{n+1} =\calU_{n+3}\setminus \{\pm 1\}, \\
    & \calT^{\frac{\beta_{n+2,n}}{2},\frac{\beta_{n+2,n}}{2}}_{n} = \calT_{n+2} \setminus \{ t_1,t_{n+2} \}.
\end{split}
\end{equation}
We notice that the set $\calU^{\frac{\beta_{n+2,n}}{2},\frac{\beta_{n+2,n}}{2}}_{n+1}$ has already been investigated in the literature (cf. \cite{Brutman78,Ibrahimoglu16}).
    
\section{The $(\beta,\gamma)$-Chebyshev (Lobatto) points are mapped equispaced points}\label{sez_mapped}
We show that the $(\beta,\gamma)$-Chebyshev (Lobatto) points can be obtained by mapping equispaced points via the so-called Kosloff Tal-Ezer (KTE) map \cite{Adcock16, Kosloff93}
\begin{equation}\label{kte}
    M_{\alpha}(x)\coloneqq \frac{\sin(\alpha\pi x/2)}{\sin(\alpha\pi /2)},\alpha\in]0,1],\;x\in\Omega.
\end{equation}
Letting
\begin{equation*}
    \mathcal{S}_n\coloneqq \bigg\{1-\frac{2j-1}{n}\bigg\}_{j=1,\dots,n},
\end{equation*}
then the Chebyshev points of first kind are
\begin{equation*}
   \mathcal{T}_n=M_1\big(\mathcal{S}_{n}\big)\coloneqq\big\{M_1(s_j)\:|\:s_j\in \mathcal{S}_{n}\big\}_{j=1,\dots,n}\,.
\end{equation*}
Moreover, if
\begin{equation*}
   \mathcal{E}_n\coloneqq \bigg\{1-\frac{2j}{n-1}\bigg\}_{j=0,\dots,n-1}
\end{equation*}
then the Chebyshev-Lobatto points are 
\begin{equation*}
    \mathcal{U}_n=M_1\big(\mathcal{E}_{n}\big)\coloneqq\big\{M_1(e_j)\:|\:e_j\in \mathcal{E}_{n}\big\}_{j=0,\dots,n}\,.
\end{equation*}
Similarly for the $(\beta,\gamma)$-Chebyshev (Lobatto) points we can prove the following result.
\begin{pro}\label{prop_mapping}
Let $\beta,\gamma\in[0,2[$, $\beta+\gamma<2$ and let $\overline{\Omega}_{\beta,\gamma}\coloneqq[-1+\beta,1-\gamma]$. Moreover, let
\begin{equation*}
   \mathcal{S}^{\beta,\gamma}_{n}\coloneqq \bigg\{1-\gamma-\frac{(2-\beta-\gamma)(2j-1)}{2n}\bigg\}_{j=1,\dots,n},
\end{equation*}
and
\begin{equation*}
    \mathcal{E}^{\beta,\gamma}_{n}\coloneqq \bigg\{1-\gamma-\frac{(2-\beta-\gamma)j}{n-1}\bigg\}_{j=0,\dots,n-1}.
\end{equation*}
Then, we have
\begin{equation*}
    \mathcal{T}^{\beta,\gamma}_{n}=M_1\big(\mathcal{S}^{\beta,\gamma}_{n}\big),\quad
    \mathcal{U}^{\beta,\gamma}_{n}=M_1\big(\mathcal{E}^{\beta,\gamma}_{n}\big).
\end{equation*}
\end{pro}
\begin{proof}
It is sufficient to observe that
\begin{equation*}
\sin\bigg(\frac{\pi}{2}\bigg(1-\gamma-\frac{(2-\beta-\gamma)(2j-1)}{2n}\bigg)\bigg) = \cos\bigg(\frac{(2-\beta-\gamma)(2j-1)\pi}{4n}+\frac{\gamma\pi}{2}\bigg)
\end{equation*}
and
\begin{equation*}
\sin\bigg(\frac{\pi}{2}\bigg(1-\gamma-\frac{(2-\beta-\gamma)j}{n-1}\bigg)\bigg) = \cos\bigg(\frac{(2-\beta-\gamma)j\pi}{2(n-1)}+\frac{\gamma\pi}{2}\bigg).    
\end{equation*}
\end{proof}
The set of $(\beta,\gamma)$-Chebyshev and $(\beta,\gamma)$-CL points are linked together as follows.
\begin{pro}\label{prop_equiv}
Let $\beta,\gamma\in[0,2[$, $\beta+\gamma<2$. Then
\begin{equation*}
    \mathcal{T}^{\beta,\gamma}_{n}=\mathcal{U}^{\beta+\frac{2-\beta-\gamma}{2n},\gamma+\frac{2-\beta-\gamma}{2n}}_{n}.
\end{equation*}
\end{pro}
\begin{proof}
In view of Proposition \ref{prop_mapping}, it is sufficient to prove the identity
\begin{equation*}
    \mathcal{S}^{\beta,\gamma}_{n}=\mathcal{E}^{\beta+\frac{2-\beta-\gamma}{2n},\gamma+\frac{2-\beta-\gamma}{2n}}_{n}.
\end{equation*}
To simplify the notation we denote as $ \rho=2-\beta-\gamma$, then for all $j=1,\dots,n$, we get
\begin{equation*}
    \begin{split}
         1-\gamma-\frac{\rho(2j-1)}{2n}&=1-\gamma-\frac{\rho}{2n}-\frac{j-1}{n-1}\bigg(\rho-\frac{\rho}{n}\bigg);\\
         -\frac{\rho(2j-1)}{2n}+\frac{\rho}{2n}&=\frac{1-j}{n-1}\rho\bigg(1-\frac{1}{n}\bigg);\\
         \frac{\rho(1-j)}{n}&=\frac{1-j}{n-1}\rho\bigg(\frac{n-1}{n}\bigg)\\
    \end{split}
\end{equation*}
and this concludes the proof.
\end{proof}
\begin{remark}
The result in Proposition \ref{prop_equiv} allows us to restrict to the set of $(\beta,\gamma)$-CL points. Moreover, we recall a remarkable property of classical Chebyshev points, that is
\begin{equation*}
    \mathcal{T}_{n}=\mathcal{U}^{\frac{1}{n},\frac{1}{n}}_{n}.
\end{equation*}
\end{remark}


\section{Lebesgue constant for $(\beta,\gamma)$-CL nodes}\label{sez_leb}
In this section we analyze the behavior of the Lebesgue constant of the $(\beta,\gamma)$-CL points on $\Omega$. 


If $\beta$ and $\gamma$ are \textit{small enough}, we can think of the $(\beta,\gamma)$-CL points as perturbed CL points and the next theorem shows that the Lebesgue constant of $(\beta,\gamma)$-CL points grows logarithmically.
\begin{thm}
Letting  $\beta,\gamma \in [0,2[$, $\beta+\gamma<2$, $\delta\coloneqq\max\{ \beta, \gamma \}$ and $\mathcal{U}^{\beta,\gamma}_{n+1}$ the set of $n+1$ $(\beta,\gamma)$-CL points, $n\in\mathbb{N}$. If
$$ \delta < \frac{4}{\pi n^2 (2+\pi \log(n+1))},$$
then
$$\Lambda(\mathcal{U}^{\beta,\gamma}_{n+1},\Omega)=\mathcal{O}(\log{n}).$$
\end{thm}

\begin{proof}
In  \cite{PiazVian18}, it has been proved that the set $\mathcal{U}_{n+1}$ under a maximal perturbation $\epsilon$ such that 
\begin{equation}\label{eq_piazzon}
\epsilon < \frac{2}{n^2(2+\pi \log(n+1))}
\end{equation}
retains the logarithmic growth of the Lebesgue constant.

In our setting, taking $\widetilde{u}_j\in\mathcal{U}^{\beta,\gamma}_{n+1}$ and $u_j\in\mathcal{U}_{n+1}$, we consider $\epsilon_j\coloneqq |\widetilde{u}_j-u_j|,\; \; j=0,\dots,n$ as the perturbation of the $j$-th point, that is

\begin{eqnarray*}
        \epsilon_j &=& \left| \cos\left( \frac{(2-\beta-\gamma)}{2n}j\pi +\frac{\gamma \pi}{2} \right) - \cos\left( \frac{j\pi}{n} \right) \right| \\
        &\leq& \left|  \frac{\beta +\gamma}{2n}j\pi -\frac{\gamma \pi}{2} \right| \leq \frac{\delta \pi}{2}\,.
\end{eqnarray*}

Then, thanks to \eqref{eq_piazzon}, if
$$ \frac{\delta \pi}{2} < \frac{2}{n^2 (2+\pi \log (n+1))},$$
the Lebesgue constant grows logarithmically, as claimed.
\end{proof}

Two particular cases deserve to be analyzed.

\subsection{Case $\beta>0,\;\gamma=0$}\label{sec_beta_leb}
We point out that we do not need to consider the case $\beta=0$ and $\gamma>0$ by virtue of the symmetric property stated in Corollary \ref{cor:symmetry}.

\begin{figure}[H]
  \centering
  \includegraphics[width=0.49\linewidth]{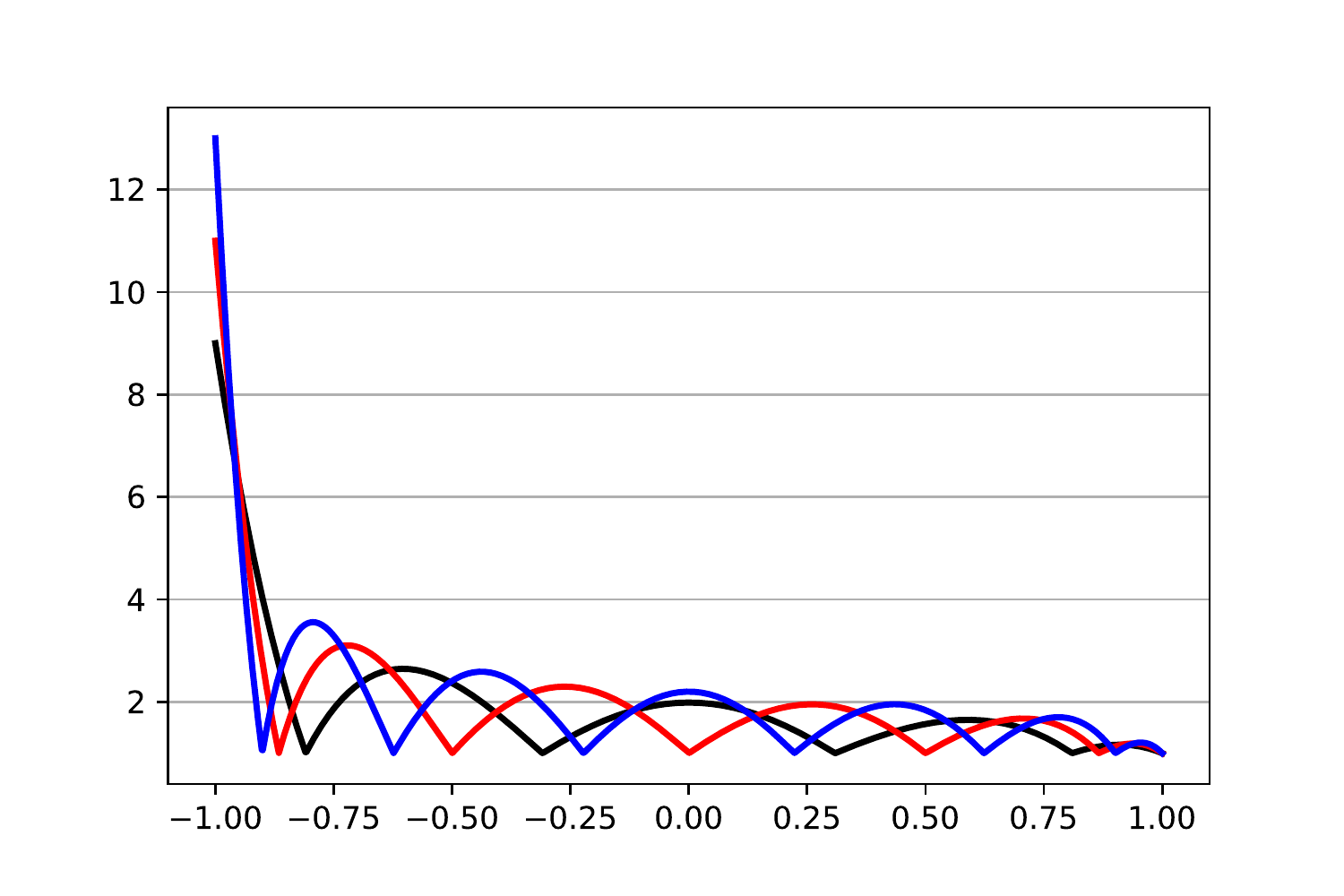}  
  \includegraphics[width=0.49\linewidth]{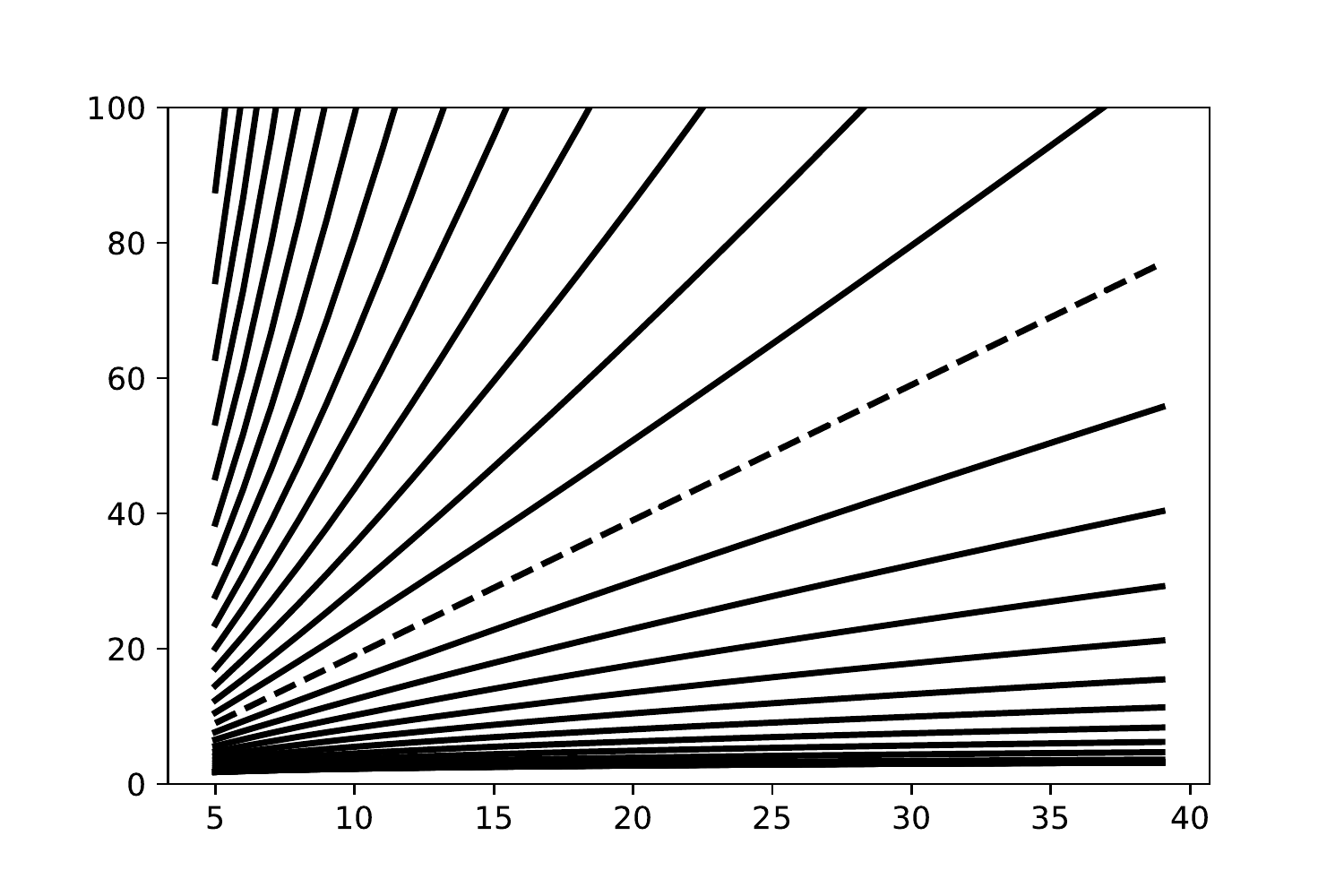}  
\caption{Left: the function $\lambda\big(\mathcal{U}^{\bar{\beta}_n,0}_{n}; \cdot \big)$ with $n=5$ (black), $n=6$ (red) and $n=7$ (blue). Right: varying $n=5,\dots,40$, the Lebesgue constant $\Lambda\big(\mathcal{U}^{\beta_{j,n},0}_{n},\Omega\big)$ with $\beta_{j,n}=j/(10n)\,,\;j=0,\dots,40$. The linear case $\beta_{20,n}=\bar{\beta}_{n}$ is displayed using a dashed line.}
\label{fig4}
\end{figure}

We start by considering the special case, $\bar{\beta}_n=2/n$, which yields to the set of points (cf. \eqref{eq_punti_rimanenti}) 
\begin{equation*}
    \mathcal{U}^{\bar{\beta}_n,0}_{n}=\mathcal{U}_{n+1}\setminus\{-1\}.
\end{equation*}
\begin{thm}\label{teor_intero}
Let $\bar{\beta}_n=2/n$, $n\in\mathbb{N}_{>0}$. Then for the associate Lebesgue function we have
\begin{equation*}
    \lambda\big(\mathcal{U}^{\bar{\beta}_n,0}_{n};-1\big)=2n-1.
\end{equation*}
\end{thm}
\begin{proof}
Let $\mathcal{U}_{n}=\{u_j\}_{j=0,\dots,n-1}$.
Then, for $i=0,\dots,{n-1}$ and $x\in\Omega$, we have
\begin{align*}
    \ell_i(x) &= \prod_{\substack{j=0 \\ j\neq i}}^{n-1} x-u_j \cdot \prod_{\substack{j=0 \\ j\neq i}}^{n-1} \frac{1}{ u_i-u_j} \\
    &= \prod_{\substack{j=0 \\ j\neq i}}^{n-1} x-u_j \cdot  (u_i-u_{n}) \prod_{\substack{j=0 \\ j\neq i}}^{n} \frac{1}{ u_i-u_j}  \\
    &= \prod_{\substack{j=0 \\ j\neq i}}^{n-1} x-u_j \cdot (u_i+1) \frac{2^{{n-1}}}{n} (-1)^i  \sigma_i,
\end{align*} 
with $\sigma_i=1/2$ if $i=0$ and $\sigma_i=1$ otherwise (see e.g. \cite[p. 37]{TrefeATAP}). Therefore,
\begin{align*}
    \ell_i(-1) &= \prod_{\substack{j=0 \\ j\neq i}}^{n-1} (-1-u_j) \cdot \frac{2^{{n-1}}}{n} (-1)^i (u_i+1) \sigma_i \\
    &= \prod_{j=0}^{n-1} (-1-u_j)\cdot\frac{1}{-1-u_i} \cdot \frac{2^{{n-1}}}{n} (-1)^i (u_i+1) \sigma_i.
\end{align*} 
We notice that 
$$\prod_{j=0}^{n-1} (-1-u_j)= \prod_{j=0}^{n-1} (u_n-u_j)$$ 
which is the inverse of the $(n+1)$-barycentric weight related to the Lagrange interpolant at the CL nodes $\mathcal{U}_{n+1}$. Hence,
$$\prod_{j=0}^{n-1} (-1-u_j) = 2^{-n+2}\,n\,(-1)^{n}.$$
This leads to
$$ |\ell_i(-1)| =  \left| \frac{1}{-1-u_i} \cdot 2^{-n+2}\,n\,(-1)^{n} \cdot \frac{2^{{n-1}}}{n} (-1)^i (u_i+1) \sigma_i \right| = 2\,\sigma_i. $$
Finally, 
$$\lambda\left(\mathcal{U}^{\bar{\beta}_n,0}_{n};-1\right) = \sum_{i=0}^{n-1} |\ell_i(-1)| = \sum_{i=0}^{n-1} 2 \, \sigma_i = 2(n-1)+1. $$
\end{proof}
\begin{conj}\label{conjetta}
In view of Theorem \ref{teor_intero}, we claim that the maximum of the Lebesgue function is attained in $x=-1$, that is
$$\Lambda\big(\mathcal{U}^{\bar{\beta}_n,0}_{n},\Omega\big)=\lambda\big(\mathcal{U}^{\bar{\beta}_n,0}_{n};-1\big)=2n-1.$$
\end{conj}
The statements in Theorem \ref{teor_intero} and Conjecture \ref{conjetta} are displayed in Figure \ref{fig4} for some values of $n$.

As supported by extensive numerical tests, the Lebesgue constant $\Lambda\big(\mathcal{U}^{\beta,0}_{n},\Omega\big)$ passes from a logarithmic to a linear growth with $n$ by increasing the value of $\beta$ from $\beta=0$ to $\beta=\bar{\beta}_n$. Moreover, as $\beta>\bar{\beta}_n$ gets larger, the growth becomes exponential. We show this behavior in Figure \ref{fig4} (right).

\begin{remark} Let $n\in\mathbb{N}$ be fixed. From numerical experiments we notice that there exists a $\beta^{\star}_n\in[0,1/n[$ such that for all $\beta\in[0,2[,\;\beta\neq\beta^{\star}_n,$
$$\Lambda\big(\mathcal{U}^{\beta^{\star}_n,0}_{n},\Omega\big)<\Lambda\big(\mathcal{U}^{\beta,0}_{n},\Omega\big).$$
Furthermore, $\Lambda\big(\mathcal{U}^{\beta,0}_{n},\Omega\big)$ is monotonically decreasing for $\beta\in[0,\beta^{\star}_n[$ and increasing for $\beta\in]\beta^{\star}_n,2[$ (see Figure \ref{fig5}). For some values of $\beta$, the growth of $\Lambda\big(\mathcal{U}^{\beta,0}_{n},\Omega\big)$ is slower than the growth of the Lebesgue constant related to the classical CL points.

\begin{figure}[H]
  \centering
  \includegraphics[width=0.49\linewidth]{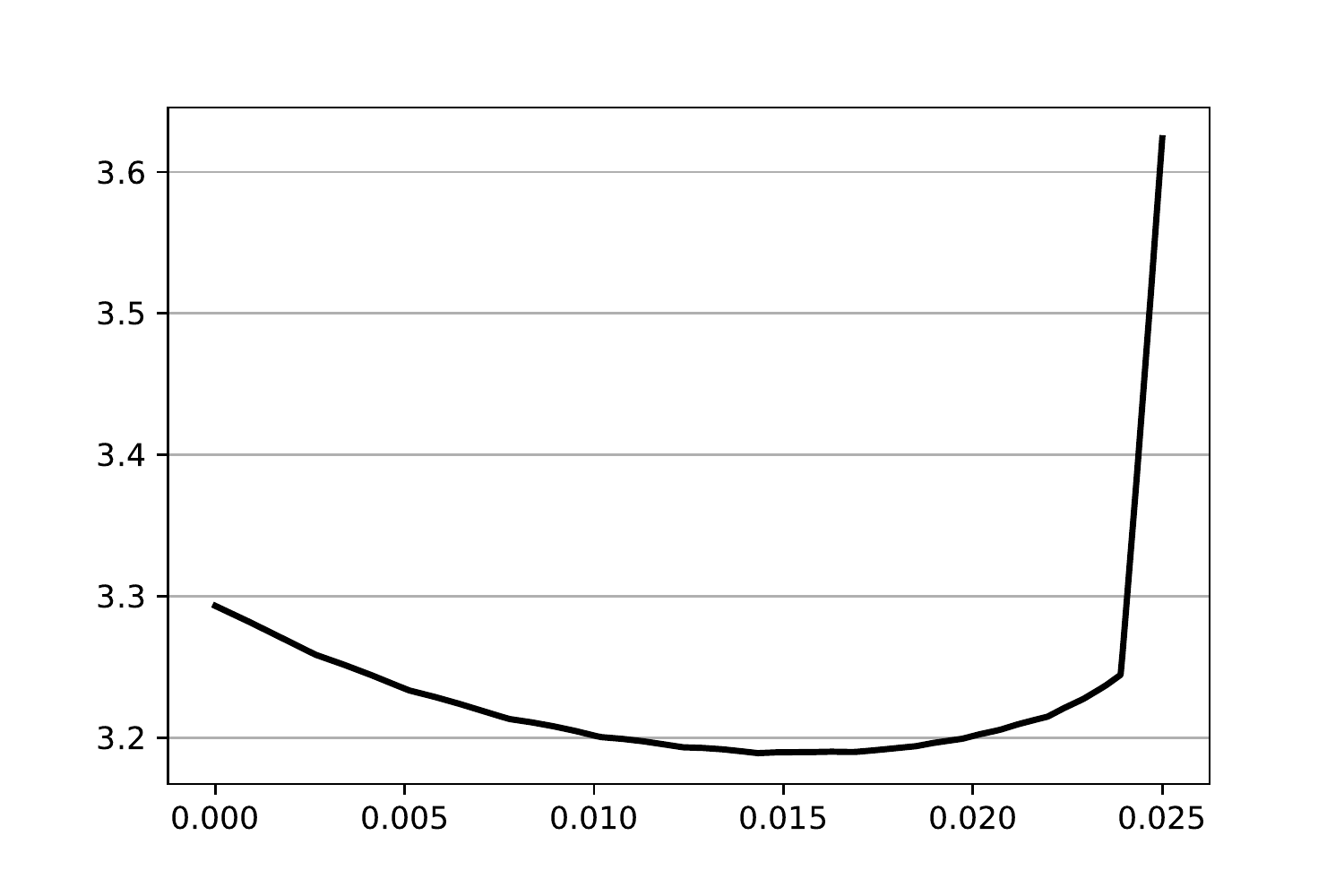}  
\includegraphics[width=0.49\linewidth]{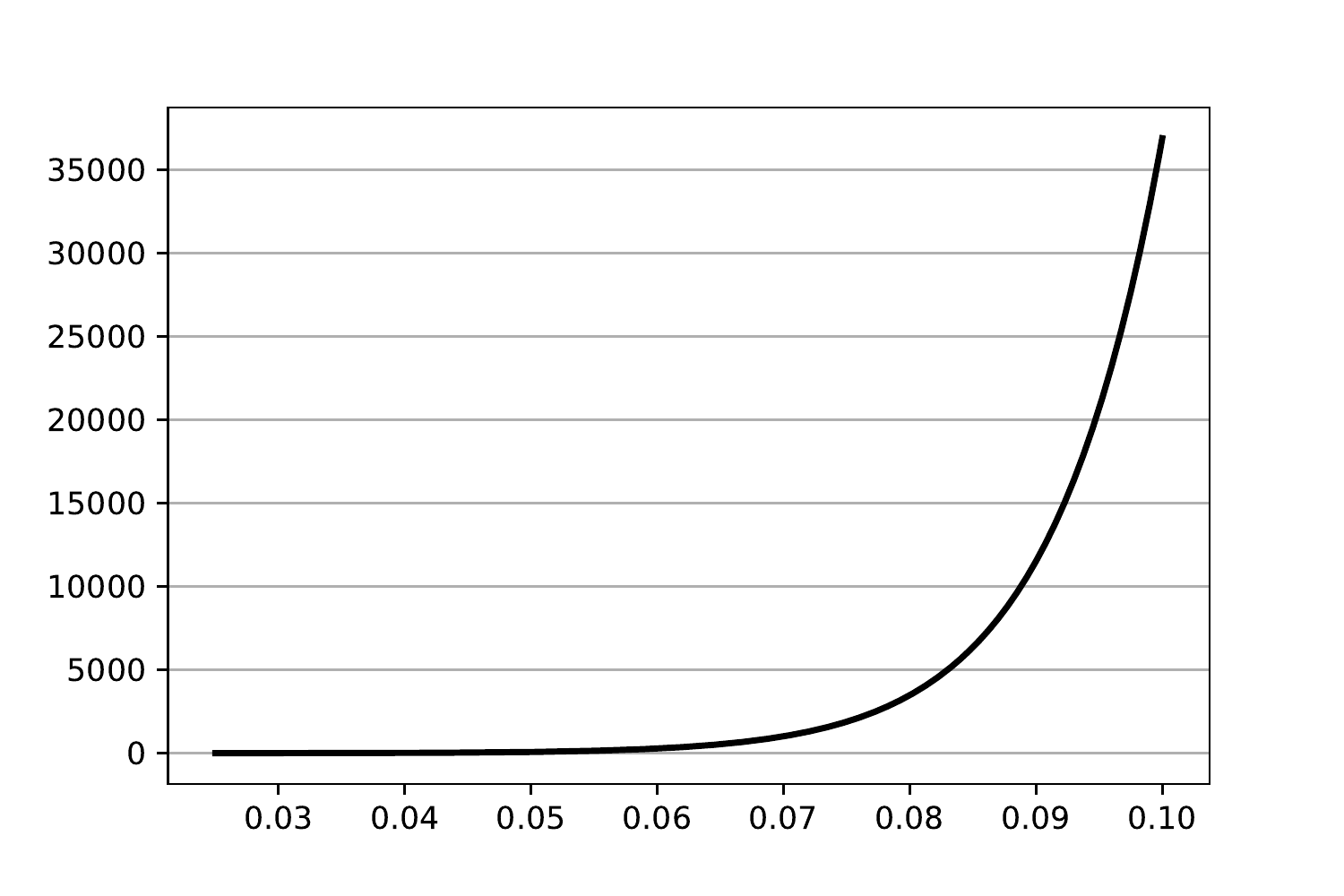}  
\caption{Fixed $n=40$, the Lebesgue constant $\Lambda\big(\mathcal{U}^{\beta,0}_{n},\Omega\big)$ varying $\beta\in[0,1/n]$ (left) and $\beta\in[1/n,4/n]$ (right).}
\label{fig5}
\end{figure}   
\end{remark}

\subsection{Case $\beta=\gamma$}
Let now $\bar{\delta}_n=2/(n+1)$. Recalling \eqref{eq_puntirimasti2}, we have that (see e.g. \cite{Brutman97})
\begin{equation*}
\Lambda\big(\mathcal{U}^{\bar{\delta}_n,\bar{\delta}_n}_{n},\Omega\big)=n.
\end{equation*}

\begin{figure}[H]
  \centering
  \includegraphics[width=0.49\linewidth]{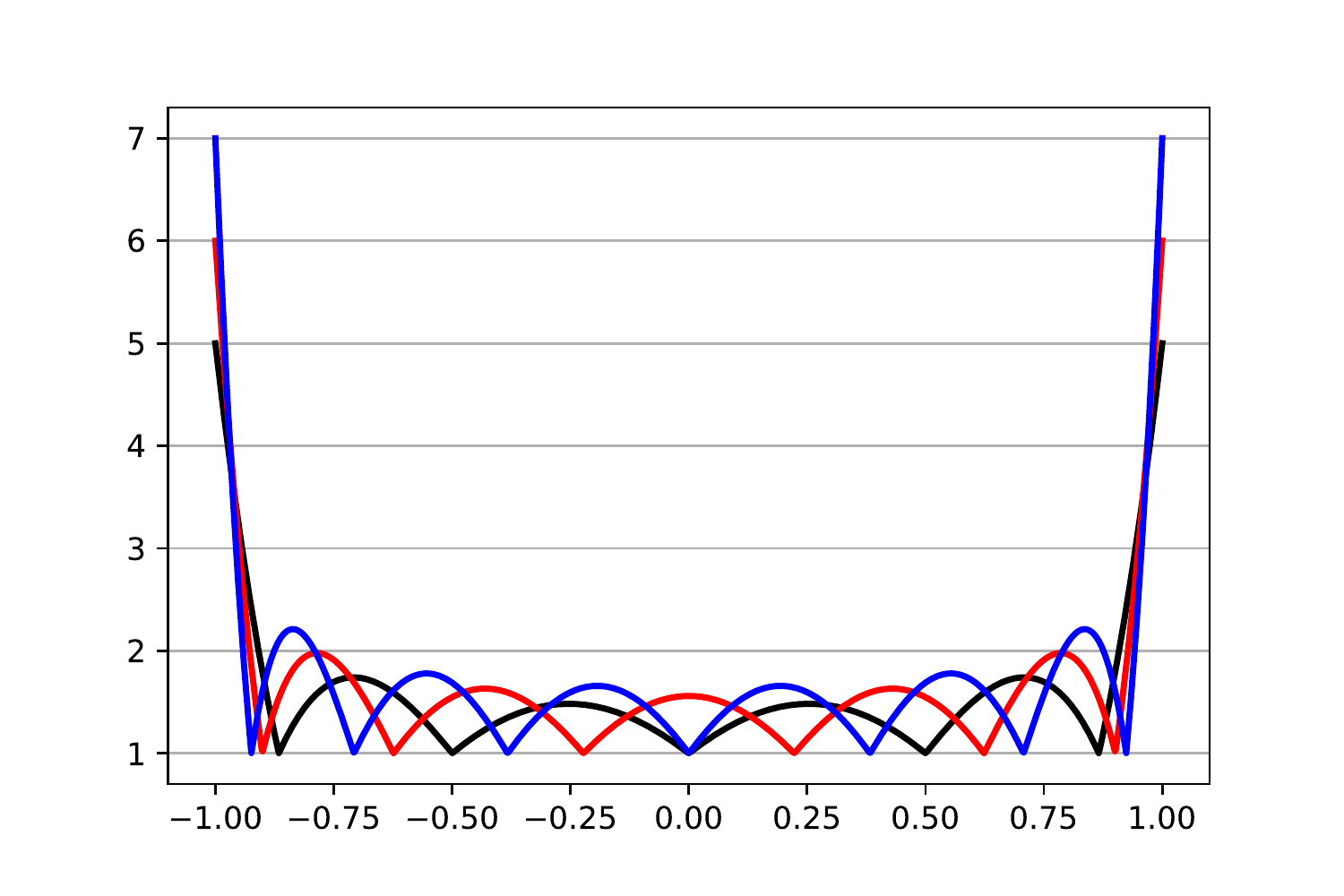}  
  \includegraphics[width=0.49\linewidth]{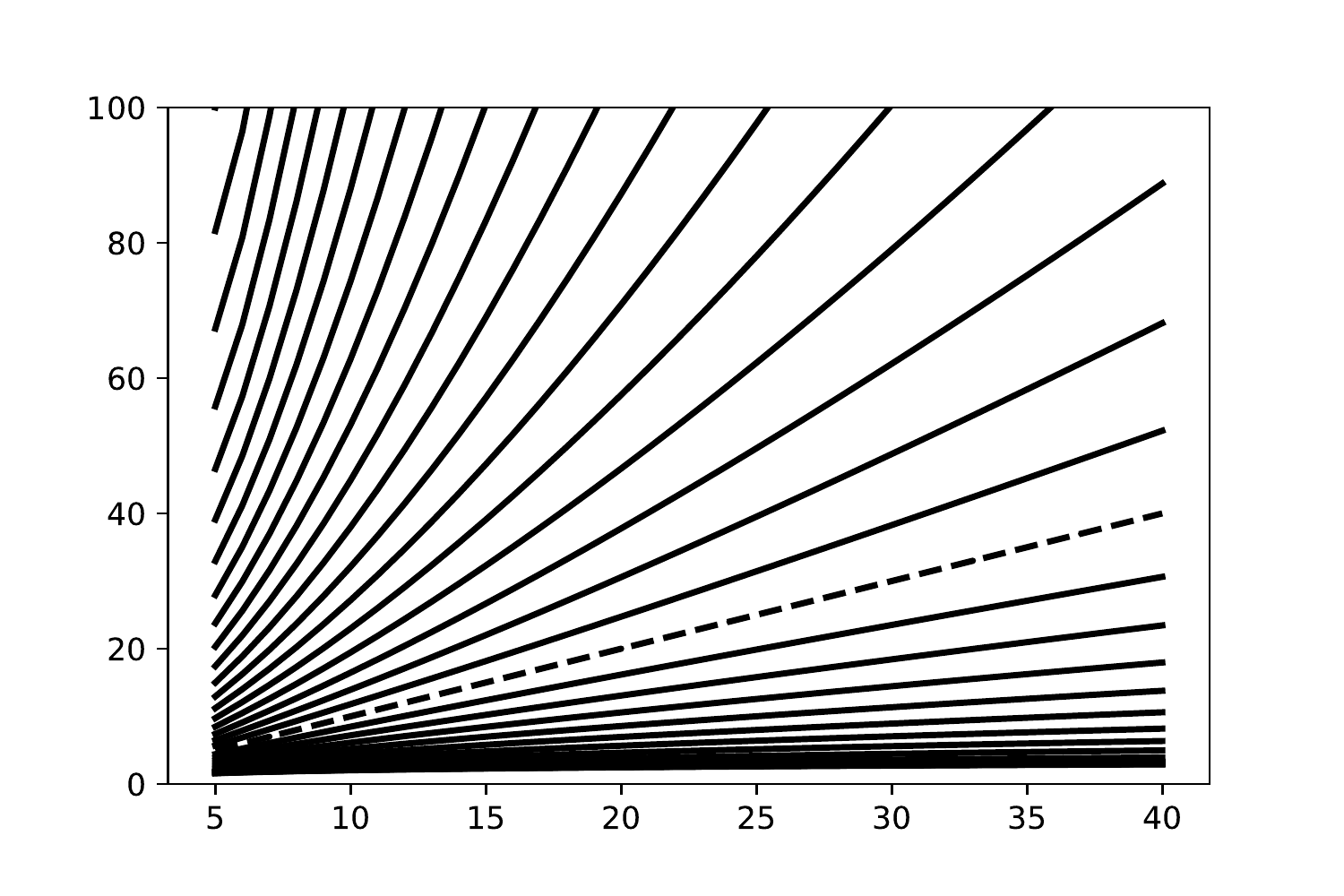}  
\caption{Left: the function $\lambda\big(\mathcal{U}^{\bar{\delta}_n,\bar{\delta}_n}_{n};\cdot\big)$ with $n=5$ (black), $n=6$ (red) and $n=7$ (blue). Right: varying $n=5,\dots,40$, the Lebesgue constant $\Lambda\big(\mathcal{U}^{\delta_{j,n},\delta_{j,n}}_{n},\Omega\big)$ with $\delta_{j,n}=j/(10(n+1))$, $j=0,\dots,40$. The linear case $\delta_{20,n}=\bar{\delta}_{n}$ is displayed using a dashed line.}
\label{fig6}
\end{figure}

As we show in Figure \ref{fig6}, considerations similar to those in Section \ref{sec_beta_leb} can be drawn, with $\bar{\delta}_n$ playing the role of $\bar{\beta}_n$.

\begin{remark}
Let $n\in\mathbb{N}$ be fixed. From numerical experiments we notice that there exists $\delta^{\star}_n\in[0,1/(n+1)[$ such that $\forall\delta\in[0,1[,\;\delta\neq\delta^{\star}_n,$
$$\Lambda\big(\mathcal{U}^{\delta^{\star}_n,\delta^{\star}_n}_{n},\Omega\big)<\Lambda\big(\mathcal{U}^{\delta,\delta}_{n},\Omega\big).$$
Moreover, $\Lambda\big(\mathcal{U}^{\delta,\delta}_{n},\Omega\big)$ is monotonically decreasing for $\delta\in[0,\delta^{\star}_n[$ and increasing for $\delta\in]\delta^{\star}_n,1[$. In Figure \ref{fig7} we plot the Lebesgue constant for different values of $\delta$. The behavior is slightly different from the results shown in Figure \ref{fig6}, in fact the minimum is achieved just before the blowing up.

\begin{figure}[H]
  \centering
  \includegraphics[width=0.49\linewidth]{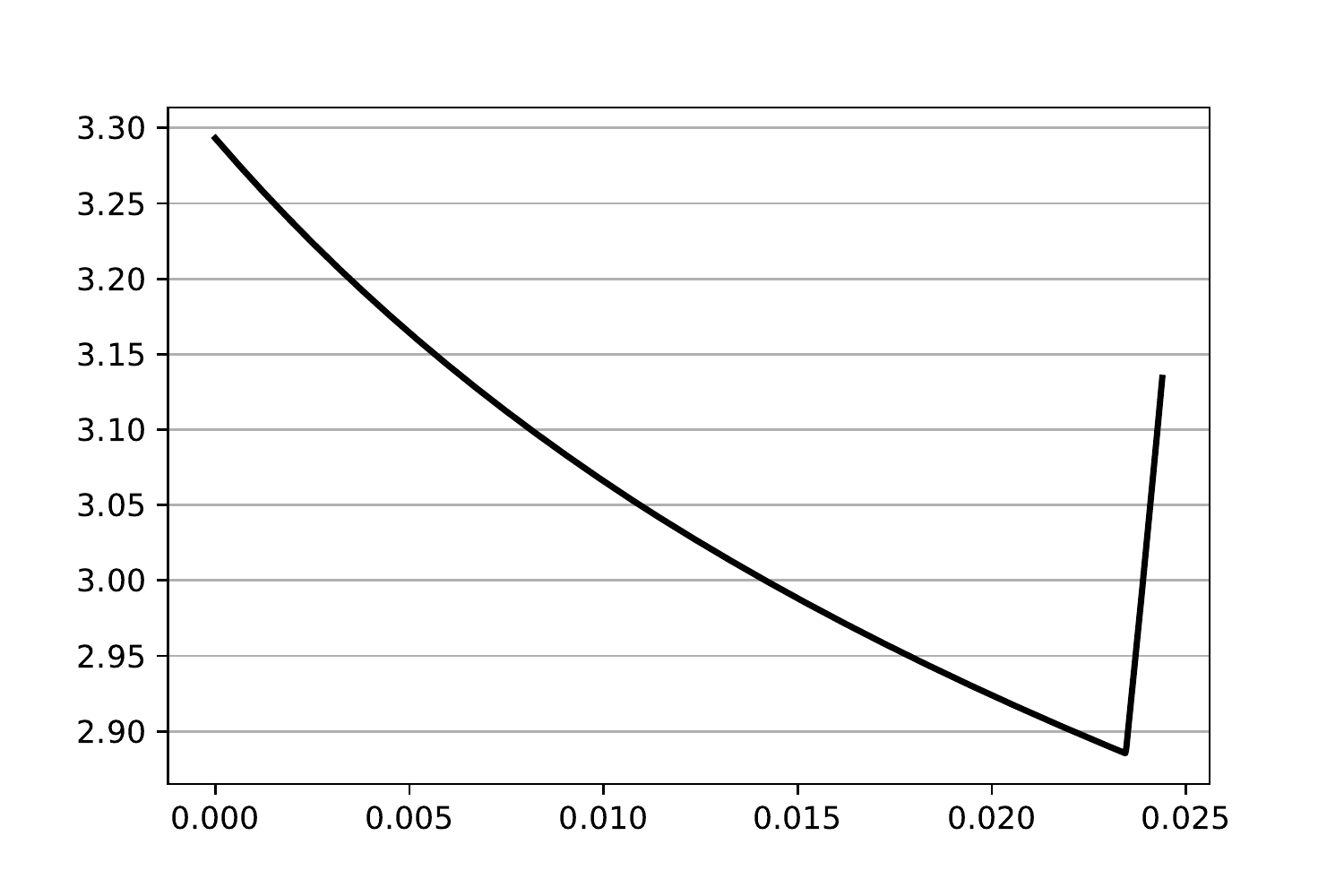}  
\includegraphics[width=0.49\linewidth]{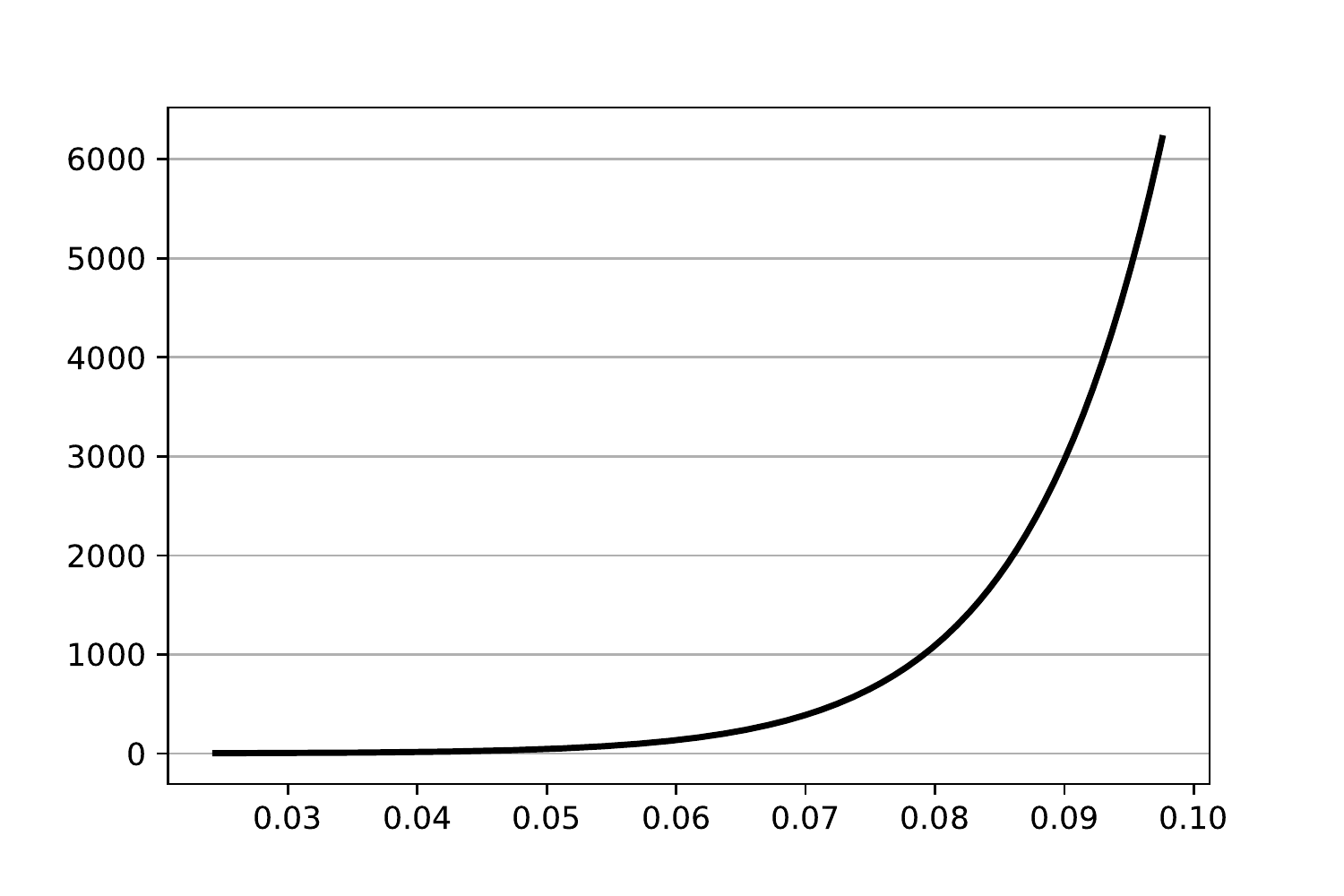}  
\caption{Fixed $n=40$, the Lebesgue constant $\Lambda\big(\mathcal{U}^{\delta,\delta}_{n},\Omega\big)$ varying $\delta\in[0,1/(n+1)]$ (left) and $\delta\in[1/(n+1),4/(n+1)]$ (right).}
\label{fig7}
\end{figure}

\end{remark}

\section{Conclusions}\label{sez_conclusions}

In this work, we introduced $(\beta,\gamma)$-Chebyshev functions and points, which can be considered a generalization of classical Chebyshev polynomials and points. In particular, for some choices of the parameters we showed that $(\beta,\gamma)$-Chebyshev functions are orthogonal polynomials in $\Omega_{\beta,\gamma}\subseteq [-1,1]$ for a proper weight function (see Theorem \ref{thm_orto}), thus they may be used in Gaussian quadrature formulae or via Newton-Côtes formulae similarly to what has been done in \cite{DeMarchi21b}.\\
Furthermore, we characterized $(\beta,\gamma)$-Chebyshev points as mapped equispaced points via KTE map and we analyzed their related Lebesgue constants showing that, for certain \textit{small} values of the parameters, they preserve the logarithmic growth, as for the classical CL points, providing alternative sets for stable polynomial approximation. This construction suggests a natural extension to the tensor product polynomial approximations. Moreover, for polynomial interpolation of total degree, we can obtain good interpolation nodes having quasi-optimal approximation properties like the well-known two-dimensional Padua points in $[-1,1]^2$ (see \cite{ Bos06, Caliari05}) or, in higher dimensions, the Lissajous points \cite{Erb15} .

\section{Acknowledgments}\label{sez_acknowledgment}

This research has been accomplished within the Rete ITaliana di Approssimazione (RITA) and the thematic group on Approximation Theory and Applications of the Italian Mathematical Union. We received the support of GNCS-IN$\delta$AM and were partially funded by the ASI - INAF grant  \lq\lq Artificial Intelligence for the analysis of solar FLARES data (AI-FLARES)\rq\rq and the NATIRESCO BIRD181249 project.

\bibliographystyle{siam}
\bibliography{bibliography.bib}
\end{document}